\documentclass{article}
\usepackage{spconf,amsmath,graphicx}
\usepackage{amssymb,amsthm}
\usepackage{fixltx2e}
\usepackage[T1]{fontenc}
\usepackage{url}
\usepackage{dsfont}
\usepackage{mathtools}
\usepackage{color}
\usepackage{graphicx}
\usepackage{bm}
\usepackage{tabularx}
\usepackage{tikz}
\usepackage{tkz-graph}
\usepackage{caption}
\usepackage{enumitem}
\usepackage{amsmath,amsfonts,amssymb,booktabs}
\usepackage[numbers,sort&compress]{natbib}
\usepackage[linesnumbered,ruled,noline]{algorithm2e}
\usepackage{supertabular}
\usetikzlibrary{calc}
\usetikzlibrary{patterns}
\usetikzlibrary{decorations.pathreplacing}
\newcommand{\g}{\raisebox{0.25ex}{\tiny $>$}}

\newcommand{\R}{\mathds{R}}

\newcommand{\norm}[1]{\lVert {#1} \rVert}
\newcommand{\Norm}[1]{\left\lVert {#1} \right\rVert}
\newcommand{\abs}[1]{\lvert {#1} \rvert}

\newcommand{\define}{\coloneqq}

\newcommand{\trans}[0]{\mathsf{T}}

\newcommand{\bw}{{\bf w}}
\newcommand{\bb}{{\bf b}}
\newcommand{\by}{{\bf y}}
\newcommand{\bz}{{\bf z}}
\newcommand{\bs}{{\bf s}}

\newcommand{\bx}{{\bf x}}
\newcommand{\bg}{{\bf g}}
\newcommand{\br}{{\bf r}}

\newcommand{\bH}{{\bf H}}
\newcommand{\bh}{{\bf h}}

\newcommand{\Real}{\operatornamewithlimits{Re}}
\newcommand{\Imag}{\operatornamewithlimits{Im}}

\newcommand{\vertex}{\node[vertex]}

\title{Joint Antenna Selection and Phase-only Beamforming using Mixed-Integer Nonlinear Programming}
\name{Tobias Fischer$^\ddag$, Ganapati Hegde$^\dagger$, Frederic
Matter$^*$, Marius Pesavento$^\dagger$, Marc E. Pfetsch$^*$, Andreas M. Tillmann$^\curlyvee$ \thanks{This work was supported by the EXPRESS project within the DFG priority program CoSIP (DFG-SPP 1798).} }
\address{
	$^*$ Department of Mathematics, Technische Universit{\"a}t Darmstadt, Germany\\
	$^\dagger$ Communication Systems Group, Technische Universit{\"a}t Darmstadt, Germany\\
	$^\ddag$ Fraunhofer Institute for Industrial Mathematics ITWM, Kaiserslautern, Germany \\
	$^\curlyvee$ Visual Computing Institute \& Chair of Operations Research, RWTH Aachen University, Germany
}

\begin{document}
%
\maketitle

\begin{abstract}
  In this paper, we consider the problem of joint antenna selection and
  analog beamformer design in downlink single-group multicast
  networks. Our objective is to reduce the hardware costs by
  minimizing the number of required phase shifters at the transmitter while
  fulfilling given distortion limits at the receivers. We formulate the problem as an
  $\ell_0$ minimization problem and devise a novel branch-and-cut based
  algorithm to solve the resulting mixed-integer nonlinear program
  to optimality. We also propose a suboptimal heuristic algorithm to solve the
  above problem approximately with a low computational complexity.
  Computational results illustrate that the solutions produced by the proposed
  heuristic algorithm are optimal in most cases. The results also indicate
  that the performance of the optimal methods can be significantly improved
  by initializing with the result of the suboptimal method.
 \end{abstract}
\section{Introduction}
Optimal antenna selection under various side constraints is a long-standing open research problem that is fundamental in a number of array processing applications \cite{antenna_selection_heath, antenna_sanayei, mimo_gore}. A traditional approach in optimizing the array geometry, commonly referred to as ``array thinning'', consists of eliminating sensors from a large uniform linear array structure with the objective to create a high resolution array with a reduced number of elements \cite{linear_array_keizer, linear_oliveri,Kupnik17}. Similarly, in upcoming Massive MIMO systems, dedicated radio frequency (RF) chains for each antenna element are no longer affordable \cite{massive_mimo_larsson, an_overview_lu}. In order to reduce hardware costs, fast switching networks may be used to adaptively select only a subset of inexpensive MIMO antennas to be connected to a reduced number of costly RF transceiver chains. Similarly, in hybrid beamforming networks based on a bank of fixed analog beamformers (ABFs) \cite{spatially_sparse_ayach, beamforming_bogale, hybrid_hegde, millimeter_wave_roh}, ge\-ne\-ral\-ly a large number of ABF outputs is available from which only a small subset can be sampled in the baseband, hence antenna selection in the beamspace domain is the underlying problem.
A major challenge in designing analog beamforming solutions for phased arrays stems from the restriction that analog phase shifters only permit phase but no magnitude variations of the antenna signals. Similarly, in digital beamforming using cheap power amplifiers with nonlinear output characteristics, low peak-to-average power ratio (PAPR) outputs are required to avoid signal distortions \cite{recursive_deng, papr_han}.

\begin{figure}
	\begin{center}
		\includegraphics[scale=1.0]{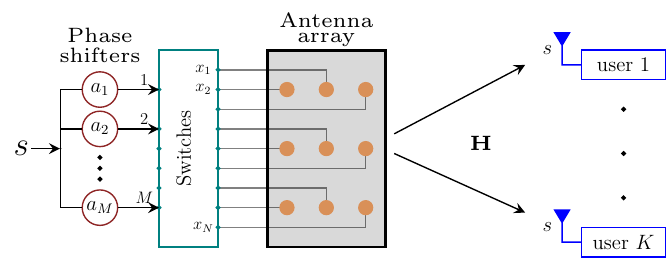}	
		\caption{Schematic diagram of system model.}
		\label{fig:system_model}
	\end{center}
\vspace{-1cm}
\end{figure}

In this paper, we consider a joint antenna selection and optimal phase-only (constant-modulus) beamforming
that can be used for array designs based on the thinning approach as well as antenna and beam selection in hybrid Massive MIMO systems.
In our approach we minimize the number of required RF phase shifters by jointly designing
the optimal phase va\-lues and assigning appropriate antenna elements to the
phase shifters, while controlling the resulting beampattern of the thinned array.
This problem is cast as an $\ell_0$ minimization task for which we
devise an efficient algorithm to solve it to (global) optimality and
demonstrate its effectiveness in numerical experiments. 
Our approach allows for an exact solution of the underlying $\ell_0$ minimization problem. We also develop a
heuristic method that computes high-quality solutions with substantially lower
computational effort. The computational results show that the solutions generated by the
proposed heuristic method are optimal in most cases, and that the heuristic
method can assist in significantly speeding up the optimal methods by providing a better initialization. 

\section{System Model}
\label{sec:system_model}
We consider a sensor array composed of $N$ antenna elements and a switching
network that selects a subset of $M$ antennas to be connected to $M$ phase shifters, where $M \ll N$,
as shown in Fig. \ref{fig:system_model}. Let $\mathbf{a} = [a_1, a_2, \ldots,
a_M]^{\trans}$ denote the ABF,
where the element $a_m \in \mathbb{C}$ is the value of the $m$th phase
shifter, for $m \in [M] \triangleq \{1,2,\ldots,M\}$. The ABF elements are
assumed to be constant modulus, i.e., $\abs{a_m} = c$ for $m \in [M]$,
where $c \in \mathbb{R}$ is a given constant. Let $\mathbf{x} \in \mathbb{C}^{N}$ be the transmit signal vector at the output of the antenna array, where its $n$th element $x_n = a_m$ if the $n$th antenna element is
connected to the $m$th phase shifter and $x_n = 0$ if the $n$th element is
not connected to any phase shifter.
There are $K$ single antenna users in the system. We assume without loss of
generality (w.l.o.g.) frequency flat channels.
Let $\mathbf{h}_k \in \mathbb{C}^{N}$ denote the channel vector between the sensor array and the $k$th user, for $k \in [K] \triangleq \{1,2,\ldots,K\}$. We define the channel matrix $\mathbf{H} = [\mathbf{h}_1, \mathbf{h}_2, \ldots, \mathbf{h}_K]$.
Let $\mathbf{s} = [s_1, s_2, \ldots, s_K]^{\trans}$ be the general form of the desired receive
signal vector at the user terminals, where the element $s_k \in \mathbb{C}$ denotes the desired
beamformer output value at the $k$th user. As a special case, in a
single-group multicast network, the desired symbol at each user is the same, i.e., $s_1=s_2=\ldots,s_K=s$ and hence $\mathbf{s} = s\mathbf{1}$. In this case, the analog beamformer can be kept constant over the coherence time of the channel, and transmit symbols can be tuned in the digital domain and be transmitted with the help of the phase shifters to serve different symbols in different time slots.

The received signal vector at the users, denoted by $\hat{\mathbf{s}} = [\hat{s}_1, \hat{s}_2, \ldots, \hat{s}_K]^{\trans}$, can be expressed as
\begin{align}
\hat{\mathbf{s}} = \mathbf{H}^{\trans}\mathbf{x} + \mathbf{n},
\end{align}
where $\mathbf{n} \in \mathcal{CN}(\mathbf{0}, \sigma^2\mathbf{I})$ represents the i.i.d. complex additive white Gaussian noise (AWGN) at the user terminals.

\section{Problem Description}
\label{sec_description}
In this section, we formulate the problem to optimally design the ABF $\mathbf{a}$ and assign antenna elements to phase shifters. Our objective is to minimize the number of required active antennas such that
the root-mean-square error between the desired receive signal vector $\mathbf{s}$ and the actual beamformer output
$\hat{\mathbf{s}}$ is at most $\sqrt{\delta}$. The problem can be mathematically
expressed as%
\begin{subequations}\label{eq:AntennaSelection}
  \begin{align}
    \underset{\mathbf{x} \in \mathbb{C}^N}{\operatorname{min\ }} &\Norm{\mathbf{x}}_0 \\
    \operatorname{s.t.\ } &\Norm{\mathbf{s} - \mathbf{H}^{\trans}\mathbf{x}}_2 \leq \sqrt{\delta}, \\
                                 &\,\abs{x_n} \in \{0, c\}, \quad \forall n \in [N],
  \end{align}
\end{subequations}
where $\norm{\bx}_0 \define \abs{\{n\in [N]: x_n\neq 0\}}$ denotes
the number of nonzero entries of~$\bx$, i.e., the number of active antennas. Throughout the remainder of this paper, we assume w.l.o.g.\ that $c = 1$.

By introducing a vector of auxiliary binary variables $\bb = [b_1, b_2, \ldots, b_N]^{\trans} \in \{0, 1\}^N$, Problem
\eqref{eq:AntennaSelection} can be reformulated equivalently with a linear
objective function as
\begin{subequations}\label{eq:BinaryReformulation}
  \begin{align}
    \underset{\bx \in \mathbb{C}^N, \bb}{\min} &\ \sum_{n = 1}^N b_n \\
    \text{s.t.}\ \quad &\norm{\mathbf{s} - \mathbf{H}^{\trans}\mathbf{x}}_2 \leq \sqrt{\delta}, \\
                         & \abs{x_n}^2 = b_n, && \forall\, n \in [N],\\
                         & b_n \in \{0, 1\}, && \forall\, n \in [N].
  \end{align}
\end{subequations}
Let $w_n \define \Real[x_n]$ and $z_n \define \Imag[x_n]$ denote the real and
ima\-ginary part of $x_n$, respectively, and let $\mathbf{w} = [w_1,w_2,
\ldots, w_N]^{\trans}$ and $\mathbf{z} = [z_1, z_2, \ldots,
z_N]^{\trans}$. Then Problem~\eqref{eq:BinaryReformulation} can be written
equivalently with real-valued variables as
\begin{subequations}
  \label{problem:real_domain_reformulation}
  \begin{align}
    \multicolumn{2}{@{}l}{$\displaystyle\underset{\bw, \bz \in \R^N, \bb}{\min}\ \sum_{n=1}^N b_n$} \\
    \text{s.t.\ } &\sum_{k=1}^K \Big( \Real[s_k] -  \left(\Real[\mathbf{h}_k]^{\trans}\mathbf{w}  - \Imag[\mathbf{h}_k]^{\trans}\mathbf{z}\right)    \Big)^2 \nonumber \\
    &+\Big( \Imag[s_k] -  \left(\Real[\mathbf{h}_k]^{\trans}\mathbf{z} + \Imag[\mathbf{h}_k]^{\trans}\mathbf{w} \right) \Big)^2 \leq {\delta},\quad \label{constraint:real_error_bound}\\
    & w_n^2 + z_n^2 \leq b_n, \quad \forall\, n \in [N], \label{eq:constraint:modulus_upper_bound}\\
    & w_n^2 + z_n^2 \geq b_n, \quad \forall\, n \in [N], \label{eq:constraint:modulus_lower_bound}\\
    & b_n \in \{0, 1\}, \,\qquad \forall\, n \in [N].
\end{align}
\end{subequations}
Note that we split the modulus constraints $\abs{x_n}^2 = w_n^2 + z_n^2 = b_n$, $n \in
[N]$, into the two inequalities~\eqref{eq:constraint:modulus_upper_bound}
and~\eqref{eq:constraint:modulus_lower_bound} since they will be handled
with different techniques.

\section{Proposed Algorithm}

Problem~\eqref{problem:real_domain_reformulation} is a (nonconvex)
mixed-integer nonlinear program (MINLP) involving integral
variables that are, in fact, binary. Moreover, it contains the
convex nonlinear (quadratic)
constraints~\eqref{constraint:real_error_bound}
and~\eqref{eq:constraint:modulus_upper_bound}, which can be rewritten as
second order cone (SOC) constraints. The nonconvexity arises due to the
quadratic constraints~\eqref{eq:constraint:modulus_lower_bound}.

MINLPs can be solved by the so-called spatial branching method, see, for
example, \cite{VigerskeGleixner2016} for a description. This method employs
a general branch-and-bound approach, which branches on integral and
continuous variables. In each node of the branch-and-bound tree, a
relaxation is solved, which might be strengthened using gradient cuts for
convex constraints. For an integral (binary) variable with a fractional
solution value, one generates two branches in which the variable is fixed
to~0 and~1, respectively. In order to guarantee satisfaction of a possibly
violated nonlinear constraint, one also creates branches on a continuous
variable by subdividing its feasible region into two parts. The reduced
regions allow to strengthen further variable bounds via so-called domain
propagation. This in turn allows to strengthen the relaxation. Under
appropriate assumptions, the method is guaranteed to converge to a global
optimum and terminates in finite time if one considers so-called
$\varepsilon$-$\delta$-feasibility, see, e.g., \cite{horst1996global}.

The general procedure can be enhanced by exploiting particular problem
structure. In the following, we show how this can be done for
Problem~\eqref{problem:real_domain_reformulation}. We describe a
customized method to propagate domains and to branch on continuous
variables using modulus constraints.
Moreover, in Section~\ref{section:heuristic_methods}, we describe a greedy heuristic method to produce upper bounds for
Problem~\eqref{problem:real_domain_reformulation}. Such bounds allow to
prune nodes in the search tree if the value of the relaxation exceeds the upper
bound. Implementations of these methods can be incorporated into an MINLP
software framework like the non-commercial solver SCIP~\cite{SCIP,
  SCIP4_0_17}, which makes it easy to employ the general spatial branching
method and add the particular customization.

\subsection{General Algorithmic Description}

In the following, we describe how to solve
Problem~\eqref{problem:real_domain_reformulation} to global optimality
using a branch-and-bound approach that is adapted to the special problem
structure.

In each node of the branch-and-bound tree, a linear programming (LP) relaxation of
Problem~\eqref{problem:real_domain_reformulation} is solved, in which $b_n
\in \{0, 1\}$ is relaxed to $0 \leq b_n \leq 1$, $n \in [N]$, and
constraints~\eqref{constraint:real_error_bound},
\eqref{eq:constraint:modulus_upper_bound}
and~\eqref{eq:constraint:modulus_lower_bound} are omitted. In
order to strengthen this linear relaxation, we add the following linear
inequalities:
\begin{align}
 -b_n \leq w_n & \leq b_n, & -b_n \leq z_n & \leq b_n,  \\
 w_n + z_n & \leq \sqrt{2}\, b_n, & w_n - z_n & \leq \sqrt{2}\, b_n, \nonumber \\
 - w_n + z_n & \leq \sqrt{2}\, b_n, & - w_n - z_n & \leq \sqrt{2}\, b_n; \nonumber
\end{align}
see Figure~\ref{fig_specbranching} for a visualization.

We use an LP relaxation in each node of the branch-and-bound tree since
this allows for fast warm-starting via the dual simplex algorithm, see,
e.g., \cite{schrijver1998}. In general, the solution of the linear
relaxation does not satisfy the nonlinear
constraints~\eqref{constraint:real_error_bound},
\eqref{eq:constraint:modulus_upper_bound}
and~\eqref{eq:constraint:modulus_lower_bound}, and the solution values of
$b_n$ need not be binary. Thus, for each binary variable with fractional solution
value in the LP relaxation, two sub-nodes are created, one for each
possible value of the binary variable. This yields tighter LP
relaxations in both sub-nodes.

As mentioned above, the error bounding
constraints~\eqref{constraint:real_error_bound} and the inequalities $w_n^2
+ z_n^2 \leq b_n$ are convex SOC constraints. If the current LP relaxation
solution violates an SOC constraint, it can be cut off by a (linear) gradient
cut~\cite{vigerske2013phd}, which is then added to the LP relaxation.

The only remaining inequalities that need to be enforced are the nonconvex
lower bound constraints $w_n^2 + z_n^2 \geq b_n$, for $n \in [N]$, to which
we refer as (lower) \emph{modulus constraints}. These nonconvex constraints
are harder to enforce. If the solution of the current LP relaxation does
not yet sa\-tisfy these constraints, we generate branching nodes, add linear
cuts or propagate domains of variables appearing in the violated modulus
constraint. These methods are described in the next subsection.

\subsection{Handling Modulus Constraints}
\label{sec:branching_rules}

\begin{figure}[t!]
  \centering
  \begin{tikzpicture}[cap=round,>=latex,every node/.style={scale=0.4}]
    \draw [thick,->] (-6, -1.6) -- (-6, 1.75) node (yaxis) [above] {\huge $z_n$};
    \draw [thick,->] (-7.6, 0) -- (-4.25, 0) node (xaxis) [right] {\huge $w_n$};
    \coordinate (s0) at (-6,0);
    \node[above left] at (-4,0) {};
    \draw[thick,blue] (s0)--++(90:1) {};
    \draw[thick,blue] (s0)--++(180:1) {};
    \draw[thick,blue] (s0)--++(-90:1) {};
    \draw[thick,blue] (s0)--++(0:1) {};
    \tikzstyle{vertex}=[circle, fill=black, draw, inner sep=1.0pt, minimum size=1.0pt]
    \vertex[] at (-6.707106781, 0.707106781) {};
    \vertex[] at (-5.292893219, -0.707106781) {};
    \vertex[] at (-6.707106781, -0.707106781) {};
    \vertex[] at (-5.292893219, 0.707106781) {};
    \vertex[] at (-6, 1) {};
    \vertex[] at (-6, -1) {};
    \vertex[] at (-5, 0) {};
    \vertex[] at (-7, 0) {};
    \draw[->, thin] (-4.9,1) -- (-4.9, 0.8);
    \draw[->, thin] (-4.9,-1) -- (-4.9, -0.8);
    \draw[->, thin] (-5,1.1) -- (-5.2, 1.1);
    \draw[->, thin] (-7,1.1) -- (-6.8, 1.1);
    \draw[->, thin] (-6.2,1.604) -- (-6.35,1.454);
    \draw[->, thin] (-6.2,-1.604) -- (-6.35,-1.454);
    \draw[->, thin] (-5.8,1.604) -- (-5.65,1.454);
    \draw[->, thin] (-5.8,-1.604) -- (-5.65,-1.454);
    \draw[cyan] (-7.25, 1) -- (-4.75, 1) node [black,right] {\huge $z_n \leq b_n$};
    \draw[cyan] (-7.25, -1) -- (-4.75, -1) node [black, right] {\huge $z_n \geq -b_n$};
    \draw[cyan] (-5, -1.25) -- (-5, 1.25) node [black, above] {\huge $w_n \leq b_n$};
    \draw[cyan] (-7, -1.25) -- (-7, 1.25) node [black, above] {\huge $w_n \geq -b_n$};
    \draw[cyan] (-4.336, -0.25) -- (-6.25, 1.654) node [black, above left]
    {\huge $w_n + z_n \leq \sqrt{2}\, b_n$};
    \draw[cyan] (-4.336, 0.25) -- (-6.25, -1.654) node [black, below left]
    {\huge $w_n - z_n \leq \sqrt{2}\, b_n$};
    \draw[cyan] (-7.664, -0.25) -- (-5.75, 1.654)  node [black, above right]
    {\huge $-w_n + z_n \leq \sqrt{2}\, b_n$};
    \draw[cyan] (-7.664, 0.25) -- (-5.75, -1.654)  node [black, below right]
    {\huge $-w_n - z_n \leq \sqrt{2}\, b_n$};
      \draw[thick] (-6cm,0cm) circle(1cm);
    \draw [thick, ->] (-2,-1.25) -- (-2,1.5) node (yaxis) [above] {\huge $z_n$};
    \draw [thick,->] (-3.2,0) -- (-0.5,0) node (xaxis) [right] {\huge $w_n$};
    \draw[thick] (-2cm,0cm) circle(1cm);
    \coordinate (s0) at (-2,0);
    \node[above left] at (-2,0) {};
    \draw[thick,blue] (s0)--++(90:1) {};
    \draw[thick,blue] (s0)--++(180:1) {};
    \draw[thick,blue] (s0)--++(-90:1) {};
    \draw[thick,blue] (s0)--++(0:1) {};
    \tikzstyle{vertex}=[circle, red,fill=red, draw, inner sep=2pt, minimum size=2pt]
    \vertex[label=above right:\huge{$(\hat{w}_n, \hat{z}_n)$}] at (-1.3, 0.6) {};
  \end{tikzpicture}
  \caption{Left: Linear inequalities for strengthening the relaxation.
    Right: Modulus constraint subdivision into orthants.}
  \label{fig_specbranching}
  \vspace{-0.4cm}
\end{figure}
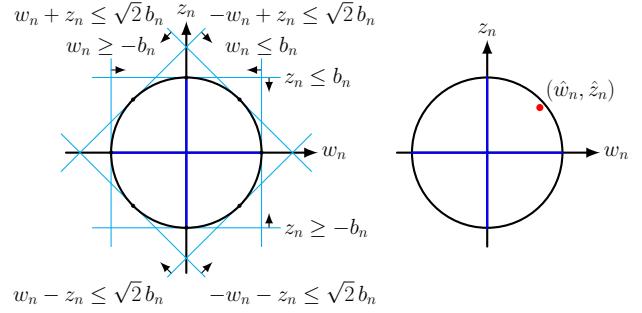

\begin{figure}[tb]
  \centering
  \newcommand{\pgfextractangle}[3]{%
    \pgfmathanglebetweenpoints{\pgfpointanchor{#2}{center}}{\pgfpointanchor{#3}{center}}
    \global\let#1\pgfmathresult
  }
  \begin{tikzpicture}[cap=round,>=latex,every node/.style={scale=0.8},scale=1.5]
    \draw [thick, ->] (0,-0.25) -- (0,2.5) node (yaxis) [above] {\large $z_n$};
    \draw [thick,->] (-0.2,0) -- (2.5,0) node (xaxis) [right] {\large $w_n$};
    \draw[thick] (2cm,0cm) arc(0:90:2cm);
    \coordinate (s0) at (0,0);
    \node[above left] at (0,0) {};
    \tikzstyle{vertex}=[circle, red,fill=red, draw, inner sep=1pt, minimum size=1pt]
    \tikzstyle{vertex2}=[circle, blue,fill=blue, draw, inner sep=1pt, minimum size=1pt]
    \coordinate (l1) at (0.3,1.97);
    \coordinate (l1x) at (0.3,0);
    \coordinate (l1y) at (0,1.97); 
    \coordinate (u1) at (1.6,1.2);
    \coordinate (u1x) at (1.6,0);
    \coordinate (u1y) at (0,1.2);
    \coordinate (l2) at (1.93,0.5);
    \coordinate (l2x) at (1.93,0);
    \coordinate (l2y) at (0,0.5);
    \coordinate (u2) at (1.05,1.7);
    \coordinate (u2x) at (1.05,0);
    \coordinate (u2y) at (0,1.7);
    \draw[thick] ($(u2y)+(-0.1,0)$) -- ($(u2y)+(0.1,0)$) node[left] at ($(u2y)+(-0.1,0)$) {$u_2$};
    \draw[thick] ($(l2y)+(-0.1,0)$) -- ($(l2y)+(0.1,0)$) node[left] at ($(l2y)+(-0.1,0)$) {$l_2$};
    \draw[thick] ($(u1x)+(0,-0.1)$) -- ($(u1x)+(0,0.1)$) node[below] at ($(u1x)+(0,-0.1)$) {$u_1$};
    \draw[thick] ($(l1x)+(0,-0.1)$) -- ($(l1x)+(0,0.1)$) node[below] at ($(l1x)+(0,-0.1)$) {$l_1$};
    \draw[thick] ($(u2x) + (0,-0.1)$) -- ($(u2x) + (0,0.1)$) node[below] at ($(u2x)+(0,-0.1)$) {$f(u_2)$};
    \draw[thick] ($(l2x) + (0,-0.1)$) -- ($(l2x) + (0,0.1)$) node[below] at ($(l2x)+(0,-0.1)$) {$f(l_2)$};
    \draw[thick] ($(u1y) + (-0.1,0)$) -- ($(u1y) + (0.1,0)$) node[left] at ($(u1y)+(-0.1,0)$) {$f(u_1)$};
    \draw[thick] ($(l1y) + (-0.1,0)$) -- ($(l1y) + (0.1,0)$) node[left] at ($(l1y)+(-0.1,0)$) {$f(l_1)$}; 

    \draw[dashed] (u2x) -- (u2);
    \draw[dashed] (u1x) -- (u1);
    \draw[dashed] (u2y) -- (u2);
    \draw[dashed] (u1y) -- (u1);

    \pgfextractangle{\angleone}{s0}{u1}
    \pgfextractangle{\angletwo}{s0}{u2}
    \draw[fill=blue, pattern=north east lines, pattern color=blue] ($(u2x)+(u1y)$) -- (u1) arc (\angleone:\angletwo:2cm);

    \vertex[label=above right:\large{$(l_1, f(l_1))$}] at (l1) {};
    \vertex[label=above right:\large{$(u_1, f(u_1))$}] at (u1) {};
    \vertex[label=above right:\large {$(f(l_2), l_2)$}] at (l2) {};
    \vertex[label=above right:\large{$(f(u_2), u_2)$}] at (u2) {};

    \draw
    [green!50!black!50,decoration={brace,amplitude=3pt},decorate]
    ($(u2x)+(0,-0.4)$) -- ($(l1x)+(0,-0.4)$) node[below, black] at (0.7,-0.45)
    {\footnotesize $l_1' = f(u_2)$};

    \draw
    [green!50!black!50,decoration={brace,amplitude=2pt},decorate]
     ($(l2x)+(0,-0.4)$) -- ($(u1x)+(0,-0.4)$) node[below, black] at (1.8,-0.45)
    {\footnotesize $u_1' = u_1$};
    \draw
    [green!50!black!50,decoration={brace,amplitude=2pt},decorate]
    ($(u2y)+(-0.65,0)$) -- ($(l1y)+(-0.65,0)$) node[left, black] at (-0.75,1.85)
    {\footnotesize $u_2' = u_2$};
    \draw
    [green!50!black!50,decoration={brace,amplitude=3pt},decorate]
    ($(l2y)+(-0.7,0)$) -- ($(u1y)+(-0.7,0)$) node[left, black] at (-0.8,0.85)
    {\footnotesize $l_2' = f(u_1)$};

  \end{tikzpicture}
  \caption{Bound propagation for the continuous variables appearing in
    modulus constraints.}
  \label{fig_propagation}
  \vspace{-0.4cm}
\end{figure}
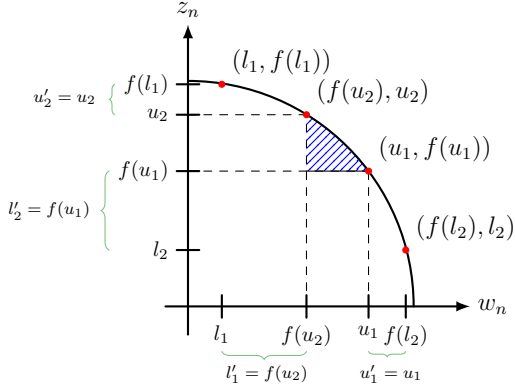

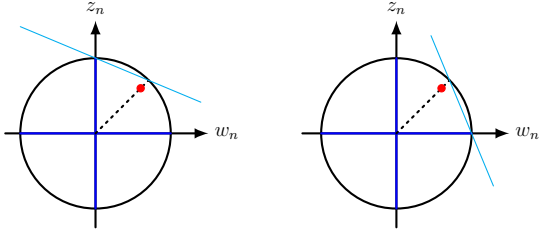
\begin{figure}[tb]
      \centering
      \begin{tikzpicture}[cap=round,>=latex,every node/.style={scale=0.65}]
        \draw [thick,->] (-4,-1.25) -- (-4,1.5) node (yaxis) [above] {\large $z_n$};
        \draw [thick,->] (-5.2,0) -- (-2.5,0) node (xaxis) [right] {\large $w_n$};
        \draw[thick] (-4cm,0cm) circle(1cm);
        \coordinate (s0) at (-4,0);
        \node[above left] at (-4,0) {};
        \draw[thick,blue] (s0)--++(90:1) {};
        \draw[thick,blue] (s0)--++(180:1) {};
        \draw[thick,blue] (s0)--++(-90:1) {};
        \draw[thick,blue] (s0)--++(0:1) {};
        \draw[thick,dotted] (s0)--++(45:1) {};
        \draw[cyan] (-5, 1.42) -- (-2.6, 0.415);
        \tikzstyle{vertex}=[circle,red, fill=red, draw, inner sep=1.5pt, minimum size=1pt]
        \vertex at (-3.4, 0.6) {};
        
        \draw [thick,->] (0,-1.25) -- (0,1.5) node (yaxis) [above] {\large $z_n$};
        \draw [thick,->] (-1.25,0) -- (1.5,0) node (xaxis) [right] {\large $w_n$};
        \draw[thick] (0cm,0cm) circle(1cm);
        \coordinate (s0) at (0,0);
        \node[above left] at (-2,0) {};
        \draw[thick,blue] (s0)--++(90:1) {};
        \draw[thick,blue] (s0)--++(180:1) {};
        \draw[thick,blue] (s0)--++(-90:1) {};
        \draw[thick,blue] (s0)--++(0:1) {};
        \draw[thick,dotted] (s0)--++(45:1) {};
        \draw[cyan] (0.46, 1.3) -- (1.29, -0.7);
        \tikzstyle{vertex}=[circle, red,fill=red, draw, inner sep=1.5pt, minimum size=1pt]
        \vertex at (0.6, 0.6) {};
      \end{tikzpicture}
      \caption{Inequalities that are added to the sub-nodes.}
      \label{fig_branchsubmodulus}
        \vspace{-0.4cm}
    \end{figure}

If the solution $(\hat{\bw}, \hat{\bz}, \hat{\bb})$ of the LP relaxation of
Problem~\eqref{problem:real_domain_reformulation} violates
$w_n^2 + z_n^2 \geq b_n$ for some $n \in [N]$, we resolve this violation
by one of the following steps:

\begin{enumerate}[leftmargin=3ex]
\item If the binary variable $\hat{b}_n$ is already fixed to
  zero, the inequality $w_n^2+z_n^2 \leq b_n$ implies that we can set
  $\hat{w}_n, \hat{z}_n$ to zero as well.
\item If the bounds of the continuous variables $w_n$ and $z_n$ are not yet
  restricted to one of the orthants w.r.t.\ $w_n\times z_n$, we create four
  branching nodes, the first with the additional constraints $w_n \geq 0$,
  $z_n \geq 0$, the second with $w_n\geq 0$, $z_n \leq 0$, the third with
  $w_n\leq 0$, $z_n \leq 0$, and the fourth with $w_n\leq 0$, $z_n \geq 0$.
  This subdivides the feasible solution set into these four orthants
  (see Figure~\ref{fig_specbranching}).
\item If the bounds of the continuous variables $w_n$ and $z_n$ are already
  restricted to one of these four orthants, we proceed as described in the following, where we
  assume w.l.o.g.\ that
  $(\hat{w}_n, \hat{z}_n, \hat{b}_n)$ is feasible for the first orthant,
  i.e., the one with $w_n \geq 0$ and $z_n \geq 0$.
  \begin{itemize}[leftmargin=3ex]
  \item[(i)] Propagation: Let $l_1 \leq w_n \leq u_1$,
    $l_2 \leq z_n \leq u_2$ denote the current lower and upper bounds of
    the variables $w_n$ and $z_n$, respectively. Compute the four points $(l_1,f(l_1))$,
    $(u_1, f(u_1))$, $(f(l_2),l_2)$ and $(f(u_2), u_2)$ on the unit circle
    that correspond to the respective lower and upper bounds of $w_n$ and $z_n$,
    where $f(x) = \sqrt{1 - x^2}$. These four points can now
    be used to strengthen the lower and upper bounds of~$w_n$ and~$z_n$. In
    order for an optimal solution $(\bw^{\star}, \bz^{\star}, \bb^{\star})$
    to fulfill the modulus constraint $w_n^2 + z_n^2 \geq b_n$, the point
    $(w_n^{\star}, z_n^{\star})$ needs to lie on or above the arc between
    the two points ($l_1',u_2')$ and $(u_1',l_2')$ if $b_n^{\star} = 1$,
    where
    \begin{align*}
      l_1' = \max \{l_1, f(u_2)\}, \quad
      u_1' = \min \{u_1, f(l_2)\}, \\
      l_2' = \max \{l_2, f(u_1)\}, \quad
      u_2' = \min \{u_2, f(l_1)\}.
    \end{align*}
    This implies that the four values $l_1'$, $u_1'$, $l_2'$ and $u_2'$ can
    now be taken as new and possibly strengthened lower and upper bounds of
    $w_n$ and $z_n$, respectively.  If the binary variable $b_n$ is not yet
    fixed to one, the lower bounds are not propagated, as $b_n$ could be
    set to zero in an optimal solution, implying $w_n = z_n = 0$ as well.
    A visualization of this propagation is given in
    Figure~\ref{fig_propagation}.
  \item[(ii)] Separation: If $\hat{w}_n + \hat{z}_n < \hat{b}_n$, add the
    cut $w_n + z_n \geq b_n$ to the LP relaxation.\footnote{
      Due to numerical reasons, we only execute (ii) if
      $\hat{w}_n^2 + \hat{z}_n^2 < 1 - \varepsilon$ with
      $\varepsilon = 10^{-5}$. Otherwise, we use standard branching rules
      for handling quadratic constraints.} Note that each
    solution in this orthant on the unit circle satisfies this inequality.
  \item[(iii)] Branching: If $\hat{w}_n + \hat{z}_n \geq \hat{b}_n$, create
    two branching nodes defined by inequalities
    $f_n\, w_n + g_n\, z_n \geq b_n$, where $f_n \in \R$ and $g_n \in \R$
    can be computed according to Figure~\ref{fig_branchsubmodulus}.
  \end{itemize}
\end{enumerate}

We prioritize the enforcement of binary variables and SOC
constraints over modulus constraints. This means that we do not use
the above methods on all modulus constraints, but only ``on demand'' in case
all the other constraints are satisfied and one particular modulus
constraint is still violated.

To select a modulus constraint for enforcing, we use a ``most infeasible''
rule: The idea is to enforce a modulus constraint $w_{\bar{n}}^2 + z_{\bar{n}}^2 \geq
b_{\bar{n}}$, $\bar{n} \in [N]$, with largest violation.  As a measure for the violation
of a modulus constraint, we use
\begin{align*}
 \rho(n) = \hat{b}_n -(\hat{w}_n^2 + \hat{z}_n^2),
\end{align*}
i.e., we choose~$\bar{n} \in [N]$ with maximal $\rho(\bar{n})$.

The whole solving procedure is summarized in Algorithm~\ref{alg1}.
Note that it is complete in the following sense: The process will terminate with a point
$(\hat{\bw}, \hat{\bz}, \hat{\bb})$ such that
$\hat{w}_n = \hat{z}_n = 0 = \hat{b}_n$ or
$1 - \varepsilon \leq \hat{w}_n^2 + \hat{z}_n^2 \leq 1$ and
$\hat{b}_n = 1$.

\begin{algorithm}[t]
\KwIn{Node of the branch-and-bound tree with current LP relaxation of the
  problem including all previously generated cuts, propagated domains and
  previously computed bounds on the objective value}
\BlankLine
obtain solution $(\hat{\bw}, \hat{\bz}, \hat{\bb})$ of LP relaxation\;
\uIf{$\hat{\bb}$ is not integral}{branch on a fractional binary variable and
  continue with another node\;}
\uElseIf{root-mean error constraint~\eqref{constraint:real_error_bound} is violated or
  $\hat{w}_n^2 + \hat{z}_n^2 > \hat {b}_n$ for some $n$}{call quadratic
  constraint handler and possibly continue with another node\;}
\uElseIf{$\hat{w}_n^2 + \hat{z}_n^2 < \hat {b}_n$}{call modulus constraint
  handler to propagate bounds or branch according to
  Section~\ref{sec:branching_rules} and continue with another
  node\;}
\Else{$(\hat{\bw}, \hat{\bz}, \hat{\bb})$ is
  optimal for the current node\;}
\caption{Node Solving Procedure within the Branch-and-Bound tree}
\label{alg1}
\end{algorithm}

\section{Heuristic Method}
\label{section:heuristic_methods}

Inspired by \cite{GESPAR} and \cite{Studer2014DemocraticR}, we propose a low-complexity suboptimal heuristic method for Problem~\eqref{eq:AntennaSelection} in the following. Let $\bh^n$ represent the
channel between the $n$th antenna element at the transmitter to all~$K$
users (i.e., the $n$th row of the channel matrix~$\bH$), and let $y_n$
denote the $n$th element of vector~$\by$. The proposed heuristic is given in Algorithm~\ref{alg2}.

We start the suboptimal algorithm with $M = 1$ active antenna
elements\footnote{We can also start with $M = M^{\text{guess}} > 1$ when a
  reasonable guess is possible or if we have any a priori knowledge about
  the minimum number of active antennas.} and increase the number of active
antenna elements by one in every iteration (in the outermost loop) until we
achieve the desired root-mean-square error bound $\sqrt{\delta}$. For each value
of $M$, we compute a large number of (max\_Iter) suboptimal solutions
$\bg_{i}$, each of them resulting from a different random
initialization of $\bx$. Note that each of the max\_Iter suboptimal
solutions $\bg_i$ can be computed in parallel to speed up the
algorithm.

As illustrated in Algorithm~\ref{alg2}, $\bx$ is updated in every iteration
(in the innermost loop) such that the root-mean-square error $e$ between
the desired and the received signal vector is decreasing. If the smallest
error $E_{i^\star}$ among all max\_Iter errors is smaller than $\sqrt{\delta}$, then the corresponding solution
$\bg_{i^\star}$ is adopted as the solution of the heuristic
method.

\begin{algorithm}[h]
  \caption{Suboptimal heuristic algorithm}
  \label{alg2}
  \KwIn{$\bH$, $\bs$, $\delta$}
  Initialize $M \leftarrow 1 \ (\text{or\ } M^{\text{guess}})$\;
  \Repeat{$E^{\star} \leq \sqrt{\delta}$ or $M > N$}{
    \For{i = $1$ to max\_Iter}{
      Randomly initialize $\bx \in \mathbb{C}^{N}$ such that
      $\norm{\bx}_{0} = M$ and $\abs{x_n} \in \{0,1\}$, $n \in [N]$\;
      Compute error $e = \norm{\bs - \bH^{\trans} \bx}_2$\;
      \For{count = $1$ to max\_Count}{
        Assign $\by \leftarrow \bx$\;
        Randomly select two integers $u$ and $v$ such that $u, v \in [N]$,
        $\abs{y_u} = 1$ and $y_v = 0$\;
        Compute the residual $\br = \bs - \bH^{\trans}
        \by + y_u (\bh^u)^{\trans}$\;
        $[{y}_u^{\star},{y}_v^{\star} ] = \underset{\bar{y}_u,
          \bar{y}_v}{\operatorname{argmin\ }} \Norm{\br -
          \bar{y}_u (\bh^u)^{\trans} + \bar{y}_v (\bh^v)^{\trans}}_2$\;
        \eIf{$\abs{{y}_u^{\star}} \geq \abs{{y}_v^{\star}}$}{
          ${y}_u \leftarrow \dfrac{{y}_u^{\star}}{\abs{{y}_u^{\star}}}$\;
        }{${y}_u \leftarrow 0$ and ${y}_v \leftarrow
          \dfrac{{y}_v^{\star}}{\abs{{y}_v^{\star}}}$\;
        }
        Compute $\hat{e} = \norm{\bs -
          \bH^{\trans}\by}_2$\;
        \If{$\hat{e} < e$ }{
          Update $\bx \leftarrow \by$ and $e \leftarrow
          \hat{e}$\;
        }
      }
      $\bg_{i} \leftarrow \bx$ and error $E_{i} \leftarrow e$\;
    }
    Compute $i^{\star} = \operatorname{argmin}_i\, E_i$\;
    $E^{\star} \leftarrow E_{i^\star}$ and $\tilde{\bx} \leftarrow
    \bg_{i^{\star}}$\;
    $M \leftarrow M + 1$\;
  }
  \Return{$\tilde{\bx}$}
\end{algorithm}

\section{Numerical Experiments}

In this section, we evaluate the performance of the proposed modulus handling based optimal method and suboptimal heuristic method in terms of the resulting hardware requirement (number of required phase shifters) and the used computation time. For the computations, we assumed a Rayleigh fading channel. The transmit
symbols are obtained by employing quadrature-phase-shift-keying (QPSK)
modulation at the transmitter with a constant magnitude $q=1.414$. We used the values $N \in
\{16, 32, 48, 64\}$ for the number of antennas, $K \in \{2, 3, 4\}$ for the
number of users and $\delta^2 \in \{0.1q, 0.2q\}$ for the error bounds.

With this setup, we solved Problem
\eqref{problem:real_domain_reformulation} with a C implementation using SCIP 4.0.1~\cite{SCIP,
  SCIP4_0_17} and CPLEX 12.7.1 as LP solver on a Linux cluster with 3.5 GHz
Intel Xeon E5-1620 Quad-Core CPUs, having 32~GB main memory and 10~MB
cache. All computations were performed single-threaded with a time
limit of one hour (3600\,s). The results are shown in
Table~\ref{Table:result}.

The table displays the solving time (in seconds) and the number of
branch-and-bound nodes for four different algorithm variants. In the
first column block, we present the results of the default version of SCIP, which
applies no special methods to handle modulus constraints -- they are
handled like general quadratic constraints. The second block shows the
results when the methods for handling modulus constraints as described in
Section~\ref{sec:branching_rules} are included in SCIP as a constraint
hand\-ler. In the third and fourth block, the results of the same two methods as
before are presented, but an initial (not necessarily
optimal) solution is computed with the suboptimal greedy heuristic method presented in
Section~\ref{section:heuristic_methods} and passed to the exact solution
method. For the number of iterations of the two inner loops, we chose
max\_Iter = max\_Count = 1000. In all four runs,
the reading times of the problem files are included in the solving times, as
are the runtimes of the suboptimal heuristic in the third and
fourth column block.
The last column block shows the sparsity $M$ of the solution
computed by the suboptimal heuristic compared to the optimal solution
computed by SCIP, as well as the solving time of the suboptimal heuristic.

The bottom part of the table presents geometric, shifted geometric
and arithmetic mean of the number of nodes and the solving
time. The shifted geometric mean of $t_1, \dots, t_n$ is defined as
\begin{align*}
\bigg( \prod_{i=1}^n (t_i + \Delta) \bigg)^{\frac{1}{n}} - \Delta,
\end{align*}
with the shift factor $\Delta = 10$ for times and $\Delta = 100$ for nodes. It
reduces the influence of easy instances on the mean values.

It turns out that the default version of SCIP already performs quite well. For $K = 2$
users, the running times are very fast even for large values of $N$. For
$K \in \{3,4\}$ users and a very small error bound
$\delta^2 = 0.1q$, the instances are much harder to solve.
From the shifted geometric mean in the bottom line, it can be seen that
adding the modulus constraint handler to SCIP results in a significantly
faster running time (about 26\,\% faster). However, the number of processed nodes does not
significantly change. The shifted geometric mean of the number of nodes that were
produced by the modulus constraint handler is 787.25 (about 24\,\%).

Executing the suboptimal heuristic and passing its solution to SCIP greatly
helps in solving the optimization problem on average, even in the default version of
SCIP. Also, the number of nodes is reduced significantly, since many nodes of
the branch-and-bound tree can be pruned. Note, however, that for the easier
problems the suboptimal heuristic consumes almost all of the solving time. Again, adding the modulus
constraint handler to SCIP speeds up the solving process (about 15\,\% and 39\,\%
speed-up compared to the default with and without initial solution, respectively), but the number of
nodes does not decrease. The shifted geometric mean of the number of nodes
produced by the modulus constraint handler is 255.43. Comparing
to the total number of nodes (500) shows
that about half of the branching nodes are used to branch on binary variables.

It is worth mentioning that the suboptimal heuristic actually returns the
optimal sparsity level in all but four instances. We observe that only for large
instances the heuristic is indeed suboptimal, but these instances cannot be solved
by SCIP within the time limit, regardless of adding the modulus constraint
handler. Interestingly, one of the instances of~Table~\ref{Table:result}
for which the heuristic computes a suboptimal solution runs into the time limit
with the default version of SCIP when this solution is passed as starting
solution. However, if the suboptimal solution is not computed beforehand, SCIP
solves this instance in roughly 700\,s. One possible explanation is the
so-called performance variability (small changes in the problem or solution
process lead to large changes in performance).

\section{Conclusion}
In this paper, we considered joint antenna selection and
design of phase-only analog beamformers. Our goal was to mi\-ni\-mize the
hardware requirement by reducing the number of phase shifters and active
antenna elements required to achieve a given maximum distortion requirement at the receivers. We
formulated the problem as an $\ell_0$ minimization program and proposed an
efficient algorithm to solve it to global op\-ti\-ma\-li\-ty. We also presented a
low-complexity suboptimal heuristic method to solve the problem approximately. The
computational results illustrate the proposed heuristic method yields
optimal solutions in most cases, with substantially reduced computational
complexity. The results also revealed that, in ge\-ne\-ral, the practical complexity of the optimal methods can be drastically reduced by initializing them with the suboptimal solutions obtained from the proposed heuristic method.

\bibliographystyle{IEEEbib}
{\small\bibliography{refs}}

\begin{thebibliography}{10}

\bibitem{antenna_selection_heath}
R.~W. Heath, S.~Sandhu, and A.~Paulraj,
\newblock ``Antenna selection for spatial multiplexing systems with linear
  receivers,''
\newblock {\em {IEEE} Commun. Letters}, vol. 5, no. 4, pp. 142--144, Apr. 2001.

\bibitem{antenna_sanayei}
S.~Sanayei and A.~Nosratinia,
\newblock ``Antenna selection in {MIMO} systems,''
\newblock {\em {IEEE} Commun. Mag.}, vol. 42, no. 10, pp. 68--73, Oct. 2004.

\bibitem{mimo_gore}
D.~A. Gore and A.~J. Paulraj,
\newblock ``{MIMO} antenna subset selection with space-time coding,''
\newblock {\em {IEEE} Trans. Signal Process.}, vol. 50, no. 10, pp. 2580--2588,
  Oct. 2002.

\bibitem{linear_array_keizer}
W.~P. M.~N. Keizer,
\newblock ``Linear array thinning using iterative {FFT} techniques,''
\newblock {\em {IEEE} Trans. Antennas and Propagation}, vol. 56, no. 8, pp.
  2757--2760, Aug. 2008.

\bibitem{linear_oliveri}
G.~Oliveri, M.~Donelli, and A.~Massa,
\newblock ``Linear array thinning exploiting almost difference sets,''
\newblock {\em {IEEE} Trans. Antennas and Propagation}, vol. 57, no. 12, pp.
  3800--3812, Dec. 2009.

\bibitem{Kupnik17}
A.~J{\"a}ger et~al.,
\newblock ``Air-coupled 40-k{H}z ultrasonic 2{D}-phased array based on a
  3{D}-printed waveguide structure,''
\newblock in {\em 2017 IEEE Int. Ultrasonics Symposium (IUS)}, Sep. 2017.

\bibitem{massive_mimo_larsson}
E.~G. Larsson, O.~Edfors, F.~Tufvesson, and T.~L. Marzetta,
\newblock ``Massive {MIMO} for next generation wireless systems,''
\newblock {\em {IEEE} Commun. Mag.}, vol. 52, no. 2, pp. 186--195, Feb. 2014.

\bibitem{an_overview_lu}
L.~Lu et~al.,
\newblock ``An overview of massive {MIMO}: Benefits and challenges,''
\newblock {\em {IEEE} J. Select. Topics in Signal Process.}, vol. 8, no. 5, pp.
  742--758, Oct. 2014.

\bibitem{spatially_sparse_ayach}
O.~E. Ayach et~al.,
\newblock ``Spatially sparse precoding in millimeter wave {MIMO} systems,''
\newblock {\em {IEEE} Trans. Wireless Commun.}, vol. 13, no. 3, pp. 1499--1513,
  Mar. 2014.

\bibitem{beamforming_bogale}
T.~E. Bogale and L.~B. Le,
\newblock ``Beamforming for multiuser massive {MIMO systems}: Digital versus
  hybrid analog-digital,''
\newblock in {\em Proc. {IEEE} Global Commun. Conf. {(GLOBECOM)}}, Austin, TX,
  USA, Dec. 2014, pp. 4066--4071.

\bibitem{hybrid_hegde}
G.~Hegde, Y.~Cheng, and M.~Pesavento,
\newblock ``Hybrid beamforming for large-scale {MIMO} systems using
  uplink-downlink duality,''
\newblock in {\em Proc. {IEEE} Int. Conf. on Acoustics, Speech and Signal
  Process. {(ICASSP)}}, Mar. 2017, pp. 3484--3488.

\bibitem{millimeter_wave_roh}
W.~Roh et~al.,
\newblock ``Millimeter-wave beamforming as an enabling technology for {5G}
  cellular communications: Theoretical feasibility and prototype results,''
\newblock {\em {IEEE} Commun. Mag.}, vol. 52, no. 2, pp. 106--113, Feb. 2014.

\bibitem{recursive_deng}
S.~K. Deng and M.~C. Lin,
\newblock ``Recursive clipping and filtering with bounded distortion for {PAPR}
  reduction,''
\newblock {\em {IEEE} Trans. Commun.}, vol. 55, no. 1, pp. 227--230, Jan. 2007.

\bibitem{papr_han}
S.~Han and J.~H. Lee,
\newblock ``{PAPR} reduction of {OFDM} signals using a reduced complexity {PTS}
  technique,''
\newblock {\em {IEEE} Signal Process. Letters}, vol. 11, no. 11, pp. 887--890,
  Nov. 2004.

\bibitem{VigerskeGleixner2016}
S.~Vigerske and A.~Gleixner,
\newblock ``{SCIP}: Global optimization of mixed-integer nonlinear programs in
  a branch-and-cut framework,''
\newblock {\em Optimization Methods and Software}, 2017,
\newblock to appear.

\bibitem{horst1996global}
R.~Horst and H.~Tuy,
\newblock {\em Global Optimization: Deterministic Approaches},
\newblock Springer, Berlin, 3 edition, 1996.

\bibitem{SCIP}
SCIP,
\newblock ``{Solving Constraint Integer Programs},'' \url{http://scip.zib.de}.

\bibitem{SCIP4_0_17}
S.~J. Maher et~al.,
\newblock ``The {SCIP} optimization suite 4.0,''
\newblock Tech. {R}ep. 17-12, ZIB, Takustr.7, 14195 Berlin, 2017.

\bibitem{schrijver1998}
A.~Schrijver,
\newblock {\em Theory of linear and integer programming},
\newblock John Wiley \& Sons, 1998.

\bibitem{vigerske2013phd}
S.~Vigerske,
\newblock {\em Decomposition in multistage stochastic programming and a
  constraint integer programming approach to mixed-integer nonlinear
  programming},
\newblock Ph.D. thesis, Humboldt-Universit{\"a}t zu Berlin, 2013.

\bibitem{GESPAR}
Y.~Shechtman, A.~Beck, and Y.~C. Eldar,
\newblock ``{GESPAR}: Efficient phase retrieval of sparse signals,''
\newblock {\em IEEE Transactions on Signal Processing}, vol. 62, no. 4, pp.
  928--938, Feb 2014.

\bibitem{Studer2014DemocraticR}
Christoph Studer, Tom Goldstein, Wotao Yin, and Richard~G. Baraniuk,
\newblock ``Democratic representations,''
\newblock CoRR abs/1401.3420, 2014.

\end{thebibliography}

\newpage
\onecolumn
\begin{center}
   \footnotesize
\setlength{\tabcolsep}{2pt}
\renewcommand{\g}{\raisebox{0.25ex}{\tiny $>$}}
\newcommand{\spc}{}
\tablehead{
\toprule
& \multicolumn{2}{@{}c@{\spc}}{default SCIP} &
\multicolumn{2}{@{\spc}c@{\spc}}{modulus handling} &
\multicolumn{2}{@{\spc}c@{\spc}}{default SCIP} &
\multicolumn{2}{@{\spc}c@{\spc}}{modulus handling} &
\multicolumn{3}{@{\spc}c@{\spc}}{suboptimal heuristic} \\
& \multicolumn{2}{@{\spc}c@{\spc}}{} &
\multicolumn{2}{@{\spc}c@{\spc}}{} &
\multicolumn{2}{@{\spc}c@{\spc}}{+ suboptimal heuristic} &
\multicolumn{2}{@{\spc}c@{\spc}}{+ suboptimal heuristic} &
\multicolumn{3}{@{\spc}c@{\spc}}{} \\
\cmidrule(l){2-3} \cmidrule(l){4-5} \cmidrule(l){6-7} \cmidrule(l){8-9} \cmidrule(){10-12}
Instance & \#nodes & time (s) & \#nodes & time (s) & \#nodes & time (s) &
\#nodes & time (s) & opt.\ sol.
& subopt.\ sol. & time (s) \\
\midrule}
\tabletail{
\midrule
\multicolumn{12}{r} \; continued on next page \\
\bottomrule
}
\tablelasttail{\bottomrule}
\tablecaption{Analysis and performance evaluation of different solution approaches/settings (SCIP 4.0.1 in optimized mode, reading time included in solving time).}
\label{Table:result}
\begin{supertabular*}{\columnwidth}{@{\extracolsep{\fill}}lrrrrrrrrrrr@{}}
$N=16, K=2, \delta^2=0.1q$, no=1 &{        469} &{        1.6}&{        501} &{        1.6} &{          7} &{        3.5} &{          9} &{        3.5} & 2 & 2 &   3.29 \\
$N=16, K=2, \delta^2=0.1q$, no=2 &{        863} &{        2.1}&{        418} &{        1.1} &{          6} &{        3.5} &{          6} &{        3.4} & 2 & 2 &   3.28 \\
$N=16, K=2, \delta^2=0.2q$, no=1 &{         11} &{        0.5}&{         11} &{        0.5} &{          1} &{        3.4} &{          1} &{        3.3} & 2 & 2 &   3.28 \\
$N=16, K=2, \delta^2=0.2q$, no=2 &{        190} &{        0.8}&{        136} &{        0.6} &{          1} &{        3.3} &{          1} &{        3.3} & 2 & 2 &   3.29 \\
$N=16, K=3, \delta^2=0.1q$, no=1 &{     3\,908} &{        9.6}&{     2\,044} &{        4.8} &{        998} &{       14.0} &{     1\,118} &{       11.4} & 4 & 4 &   9.47 \\
$N=16, K=3, \delta^2=0.1q$, no=2 &{     1\,690} &{        4.9}&{     2\,337} &{        4.6} &{        996} &{       14.0} &{        928} &{       12.1} & 4 & 4 &   9.42 \\
$N=16, K=3, \delta^2=0.2q$, no=1 &{     1\,704} &{        7.9}&{     1\,529} &{        3.5} &{         47} &{        8.8} &{         71} &{        8.7} & 3 & 3 &   7.08 \\
$N=16, K=3, \delta^2=0.2q$, no=2 &{     2\,688} &{        7.5}&{     1\,720} &{        3.3} &{         87} &{        8.4} &{         95} &{        8.3} & 3 & 3 &   7.10 \\
$N=16, K=4, \delta^2=0.1q$, no=1 &{    25\,507} &{      121.8}&{    22\,684} &{       40.7} &{     2\,374} &{       27.7} &{     2\,470} &{       20.9} & 5 & 5 &  16.39 \\
$N=16, K=4, \delta^2=0.1q$, no=2 &{     2\,853} &{       10.4}&{     2\,603} &{        6.5} &{         95} &{       14.9} &{        133} &{       14.9} & 4 & 4 &  13.04 \\
$N=16, K=4, \delta^2=0.2q$, no=1 &{     2\,533} &{       11.0}&{     1\,481} &{        3.9} &{        289} &{       15.6} &{        393} &{       15.6} & 4 & 4 &  13.06 \\
$N=16, K=4, \delta^2=0.2q$, no=2 &{    12\,676} &{       57.2}&{    11\,591} &{       20.7} &{    13\,095} &{       73.6} &{     9\,637} &{       28.2} & 5 & 5 &  16.41 \\
$N=32, K=2, \delta^2=0.1q$, no=1 &{        170} &{        7.5}&{        204} &{        9.6} &{          4} &{        4.2} &{          4} &{        4.2} & 2 & 2 &   4.04 \\
$N=32, K=2, \delta^2=0.1q$, no=2 &{     1\,937} &{        6.2}&{        714} &{        2.2} &{          2} &{        4.2} &{          2} &{        4.2} & 2 & 2 &   3.97 \\
$N=32, K=2, \delta^2=0.2q$, no=1 &{        201} &{        2.7}&{        203} &{        1.4} &{          6} &{        4.4} &{          8} &{        4.3} & 2 & 2 &   3.98 \\
$N=32, K=2, \delta^2=0.2q$, no=2 &{        101} &{        2.8}&{         41} &{        1.6} &{         17} &{        4.5} &{         11} &{        4.4} & 2 & 2 &   3.99 \\
$N=32, K=3, \delta^2=0.1q$, no=1 &{    12\,910} &{       52.1}&{     5\,582} &{       16.5} &{         77} &{       13.5} &{        134} &{       12.9} & 3 & 3 &   8.53 \\
$N=32, K=3, \delta^2=0.1q$, no=2 &{        313} &{        1.7}&{        356} &{        1.6} &{        177} &{       13.1} &{        210} &{       14.3} & 3 & 3 &   8.58 \\
$N=32, K=3, \delta^2=0.2q$, no=1 &{     3\,387} &{       12.0}&{     2\,671} &{        9.6} &{        159} &{       13.9} &{        198} &{       13.8} & 3 & 3 &   8.52 \\
$N=32, K=3, \delta^2=0.2q$, no=2 &{     1\,380} &{        7.2}&{     2\,514} &{       10.7} &{        143} &{       14.1} &{        151} &{       14.2} & 3 & 3 &   8.61 \\
$N=32, K=4, \delta^2=0.1q$, no=1 &{    96\,569} &{      559.3}&{    13\,703} &{       51.0} &{   154\,122} &{      652.4} &{    57\,158} &{      156.6} & 4 & 5 &  19.59 \\
$N=32, K=4, \delta^2=0.1q$, no=2 &{   141\,531} &{      816.5}&{    97\,762} &{      301.4} &{   147\,573} &{      648.7} &{    66\,224} &{      172.8} & 5 & 5 &  18.29 \\
$N=32, K=4, \delta^2=0.2q$, no=1 &{     8\,070} &{       31.9}&{     9\,874} &{       37.2} &{         43} &{       15.6} &{         42} &{       15.0} & 3 & 3 &  11.75 \\
$N=32, K=4, \delta^2=0.2q$, no=2 &{     4\,879} &{       20.5}&{    14\,224} &{       52.9} &{         97} &{       16.5} &{        123} &{       16.7} & 3 & 3 &  11.73 \\
$N=48, K=2, \delta^2=0.1q$, no=1 &{        496} &{        8.7}&{        315} &{        8.2} &{          1} &{        4.9} &{          1} &{        4.9} & 2 & 2 &   4.64 \\
$N=48, K=2, \delta^2=0.1q$, no=2 &{        190} &{        2.0}&{        734} &{        5.0} &{          9} &{        5.1} &{          9} &{        5.1} & 2 & 2 &   4.59 \\
$N=48, K=2, \delta^2=0.2q$, no=1 &{         15} &{        1.6}&{        447} &{        3.1} &{          7} &{        5.0} &{          7} &{        5.0} & 2 & 2 &   4.58 \\
$N=48, K=2, \delta^2=0.2q$, no=2 &{        201} &{        2.9}&{        487} &{        4.5} &{          7} &{        4.9} &{          7} &{        4.9} & 2 & 2 &   4.58 \\
$N=48, K=3, \delta^2=0.1q$, no=1 &{     2\,256} &{       18.9}&{     8\,253} &{       49.2} &{        259} &{       16.3} &{        283} &{       16.4} & 3 & 3 &   9.85 \\
$N=48, K=3, \delta^2=0.1q$, no=2 &{    39\,191} &{      236.6}&{     5\,542} &{       31.4} &{        285} &{       19.7} &{        332} &{       20.0} & 3 & 3 &   9.90 \\
$N=48, K=3, \delta^2=0.2q$, no=1 &{     1\,810} &{       21.8}&{     2\,777} &{       19.5} &{        381} &{       20.5} &{        484} &{       20.8} & 3 & 3 &   9.87 \\
$N=48, K=3, \delta^2=0.2q$, no=2 &{     3\,648} &{       27.2}&{     2\,207} &{       11.2} &{        347} &{       23.1} &{        362} &{       22.9} & 3 & 3 &   9.91 \\
$N=48, K=4, \delta^2=0.1q$, no=1 &{   104\,244} &{      992.9}&{    56\,336} &{      354.9} &{   168\,481} &{     2477.0} &{   107\,320} &{      439.4} & 4 & 5 &  22.63 \\
$N=48, K=4, \delta^2=0.1q$, no=2 &{    70\,950} &{      465.3}&{    28\,474} &{      177.2} &{     4\,895} &{       64.2} &{     6\,357} &{       62.8} & 4 & 4 &  17.99 \\
$N=48, K=4, \delta^2=0.2q$, no=1 &{     9\,764} &{       65.4}&{    29\,631} &{      170.3} &{         83} &{       24.2} &{         91} &{       24.6} & 3 & 3 &  13.53 \\
$N=48, K=4, \delta^2=0.2q$, no=2 &{    56\,580} &{      373.8}&{    67\,515} &{      348.5} &{     9\,933} &{       86.9} &{    11\,707} &{       76.0} & 4 & 4 &  18.09 \\
$N=64, K=2, \delta^2=0.1q$, no=1 &{        397} &{        4.2}&{        273} &{        3.3} &{         10} &{        5.7} &{         10} &{        5.7} & 2 & 2 &   5.18 \\
$N=64, K=2, \delta^2=0.1q$, no=2 &{        360} &{        4.0}&{        360} &{        3.9} &{         11} &{        5.7} &{         11} &{        5.7} & 2 & 2 &   5.17 \\
$N=64, K=2, \delta^2=0.2q$, no=1 &{        505} &{        4.5}&{        931} &{        7.0} &{          8} &{        5.7} &{          8} &{        5.7} & 2 & 2 &   5.19 \\
$N=64, K=2, \delta^2=0.2q$, no=2 &{        476} &{        4.8}&{        481} &{        4.8} &{         11} &{        6.0} &{         11} &{        6.0} & 2 & 2 &   5.19 \\
$N=64, K=3, \delta^2=0.1q$, no=1 &{    11\,498} &{      105.5}&{     5\,982} &{       36.5} &{        529} &{       33.0} &{        568} &{       32.6} & 3 & 3 &  11.17 \\
$N=64, K=3, \delta^2=0.1q$, no=2 &{     6\,464} &{       70.6}&{    20\,879} &{      105.7} &{        279} &{       21.9} &{        333} &{       21.9} & 3 & 3 &  11.17 \\
$N=64, K=3, \delta^2=0.2q$, no=1 &{        830} &{       12.2}&{        884} &{       11.9} &{        519} &{       34.2} &{        512} &{       34.5} & 3 & 3 &  11.24 \\
$N=64, K=3, \delta^2=0.2q$, no=2 &{     6\,626} &{       59.2}&{     9\,821} &{       57.0} &{        815} &{       34.8} &{        966} &{       34.2} & 3 & 3 &  11.22 \\
$N=64, K=4, \delta^2=0.1q$, no=1 &{   217\,176} &{     3128.0}&{   360\,779} &{     2932.9} &{    20\,299} &{      172.6} &{    24\,143} &{      152.2} & 4 & 4 &  20.38 \\
$N=64, K=4, \delta^2=0.1q$, no=2 &{    75\,975} &{      670.0}&{   329\,245} &{     2982.6} &{\g 157\,627} &{\g   3600.0} &{   510\,570} &{     3492.6} & 4 & 5 &  25.50 \\
$N=64, K=4, \delta^2=0.2q$, no=1 &{    51\,999} &{      688.0}&{    47\,384} &{      396.5} &{    10\,226} &{      126.6} &{    28\,560} &{      185.1} & 3 & 4 &  20.58 \\
$N=64, K=4, \delta^2=0.2q$, no=2 &{    82\,051} &{      748.6}&{    26\,076} &{      173.8} &{        119} &{       35.1} &{        123} &{       35.4} & 3 & 3 &  15.24 \\
\midrule
Geometric mean         &   2\,796 &     21.6&   2\,816 &     16.0&      164 &     19.6&      178 &     17.3\\
Shifted geometric mean &   3\,308 &     36.5&   3\,228 &     27.5&      470 &     26.4&      500 &     22.4\\
Arithmetic mean        &  22\,296 &    197.3&  25\,014 &    176.8&  14\,490 &    175.6&  17\,331 &    110.0\\
\end{supertabular*}
\end{center}
\twocolumn


\end{document}